\documentclass[a4paper]{amsart}
\usepackage{epsfig}
\usepackage{amsmath}
\usepackage{amsthm}
\usepackage{amssymb}
\usepackage{latexsym,amsfonts}
\usepackage{amscd}
\usepackage[all]{xy}

\usepackage{psfrag}  

\parskip\smallskipamount
\numberwithin{equation}{section}
\newcommand{\GL}{\mathrm{GL}}

\newcommand{\cO}{{\mathcal{O}}}
\newcommand{\cT}{{\mathcal{T}}}

\newcommand{\frl}{{\mathfrak l}}

\newcommand{\frn}{{\mathfrak n}}
\newcommand{\frp}{{\mathfrak p}}

\newcommand{\frgl}{{\mathfrak{gl}}}

\newtheorem{prop}{Proposition}[section]
\newtheorem{thm}[prop]{Theorem}
\newtheorem*{theorem}{Theorem}
\newtheorem{lemma}[prop]{Lemma}
\newtheorem{cor}[prop]{Corollary}
\newtheorem{ex}[prop]{Example}
\newtheorem{remark}[prop]{Remark}

\newtheorem{definition}[prop]{Definition}

\def\rk{\operatorname {rk}}

\begin{document}
\title{On the complement of the Richardson orbit}


\author[Baur]{Karin Baur}
\address{Department of Mathematics \\
ETH Z\"urich
R\"amistrasse 101 \\
CH-8092 Z\"urich \\
Switzerland
}
\email{baur@math.ethz.ch}

\author[Hille]{Lutz Hille}
\address{Mathematisches Institut \\
Fachbereich Mathematik und Informatik der Universit\"at M\"unster \\
Einsteinstrasse 62 \\
D-48149  M\"unster \\
Germany
}
\email{lutz.hille@uni-muenster.de}



\thanks{This research was supported through the programme
``Research in Pairs''
by the Mathematisches Forschungsinstitut Oberwolfach in 2009.
The second author was supported by the DFG
priority program SPP 1388 representation theory.}

\begin{abstract}
We consider parabolic subgroups of a general algebraic group
over an algebraically closed field $k$ whose Levi part has
exactly $t$ factors. By a classical theorem of
Richardson, the nilradical of a parabolic subgroup $P$
has an open dense $P$-orbit.  In the complement to this
dense orbit, there are
infinitely many orbits as soon as the number $t$ of factors
in the Levi part is $\ge 6$.
In this paper, we describe the irreducible components of
the complement. In particular, we show that there are at most
$t-1$ irreducible components. 
We are also able
to determine their codimensions.
\end{abstract}

\maketitle

\tableofcontents

%
\section{Introduction and notations}
%

Let $P$ be a parabolic subgroup of a reductive algebraic group $G$ over 
an algebraically closed field $k$. Let $\frp$ be its Lie algebra and let 
$\frp = \frl \oplus \frn$ be the Levi decomposition of $\frp$, i.e. $\frn$ is 
the nilpotent radical of $\frp$. A classical result of Richardson~\cite{ri} says that 
$P$ has an open dense orbit in the nilradical. We will call this $P$-orbit the 
{\em Richardson orbit for $P$}. 
However, in general there are infinitely many $P$-orbits in $\frn$. 

For classical $G$, the 
cases where there are finitely many $P$-orbits 
in $\frn$ have been classified in~\cite{hr1}. 
Also, the $P$-action on the derived 
Lie algebras of $\frn$ have been studied in a series of papers, and the cases 
with finitely many orbits have been classified, cf.~\cite{bh1},~\cite{bh2},~\cite{bh3},
~\cite{bhr}. 

If $G$ is a general linear group, $G=\GL_n$, then the parabolic subgroup 
$P$ can be described by the lengths of the blocks in the Levi factor: 
Write $P=L N$ where $L$ is a Levi factor and $N$ is the unipotent radical 
of $P$. Then we can assume that $L$ consists of matrices which have non-zero 
entries in square blocks on the diagonal. Similarly, the Levi factor $\frl$ of $\frp$ consists 
of the $n\times n$-matrices with non-zero entries lying in squares of size $d_i\times d_i$ 
($i=1,\dots,t$) 
on the diagonal and $\frn$ are the matrices which only have non-zero entries 
above and to the right of these square blocks. 

Let $t$ be the number of such 
blocks and $d_1,\dots,d_t$ the lengths of them, $\sum d_i=n$ 
(with $d_i>0$ for all $i$). So $d$ is a composition of $n$. 
We will 
call such a $d=(d_1,\dots,d_t)$ a {\em dimension vector}. 
We write $P(d)$ for the corresponding parabolic subgroup 
and $\frn(d)$ for the nilpotent radical of $P(d)$, the Richardson orbit of $P(d)$ is 
denoted by $\cO(d)$. Its partition will be $\lambda(d)$.
Once $d$ is fixed, we will often just use $P$, $\frn$
and $\lambda$ if there is no ambiguity. 
Recall that the nilpotent $\GL_n$-orbits are parametrised by partitions of 
$n$. We will use $C(\mu)$ to denote the nilpotent $\GL_n$-orbit for 
the partition $\mu$ ($\mu$ a partition of $n$). 
And we will usually denote $P$-orbits in $\frn$ by a calligraphic O, i.e. we will 
write $\cO$ or $\cO(\mu)$ if $\mu$ is the partition of the
nilpotency class of the $P$-orbit. 

Now, the nilradical $\frn$ is a disjoint union of the intersections $\frn\cap C(\mu)$ 
of the nilradical with all nilpotent $\GL_n$-orbits. 
By Richardsons result, $\frn\cap C(\lambda)=\cO(\lambda)$ is a single $P$-orbit.
In particular, the Richardson orbit consists exactly of the elements of the nilpotency 
class $\lambda$. 
However, for $\mu\le \lambda$, the 
intersection $\frn\cap C(\mu)$ might be reducible 
(cf. Proposition~\ref{prop:comps}). 

In the case where $\frn$ is the nilradical of a Borel subalgebra of the Lie algebra 
of a simple 
algebraic group $G$, Spaltenstein has first studied 
the varieties $\frn\cap (G\cdot e)$ for 
$G\cdot e$ a nilpotent orbit under the adjoint action (\cite{sp}). 
In \cite{ghr}, the authors study the action of a Borel subgroup $B$ of a simple 
algebraic group on the closure $\frn\cap C(\mu)$ for the subregular nilpotency 
class $C(\mu)$ and characterize the cases 
where $B$ has only finitely many orbits under the adjoint action.


The main goal of this article is to describe the irreducible components 
of the complement $Z:=\frn\setminus \cO(d)$
of the Richardson orbit in $\frn$. 
They occur in intersections $\frn\cap C(\mu)$ 
for certain partitions 
$\mu=\mu(i,j)\le\lambda$. 

We have two descriptions of the irreducible components of $Z$. 
On one hand, we give rank conditions on the matrices of $\frn$, 
on the other hand, we use tableaux $T(i,j)$ for certain 
$(i,j)$ with $1\le i<j\le t$ and 
associate irreducible components $\frn(T(i,j))$ 
of the intersections $\frn\cap C(\mu(i,j))$ 
to them. 
Before we can 
state the two results we now introduce the necessary notation. 

\begin{figure}[h]
\begin{center}
  \psfragscanon
   \psfrag{A}{$A$}
   \psfrag{12}{$_{12}$}
   \psfrag{13}{$_{13}$}
   \psfrag{11}{$_{11}$}
   \psfrag{21}{$_{21}$}
   \psfrag{22}{$_{22}$}
   \psfrag{23}{$_{23}$}
   \psfrag{33}{$_{33}$}
   \psfrag{32}{$_{32}$}
   \psfrag{31}{$_{31}$}
   \psfrag{1}{$_{1}$}
   \psfrag{2}{$_{2}$}
   \psfrag{3}{$_{3}$}
   \psfrag{d}{$_{d}$}
\includegraphics[scale=.5]{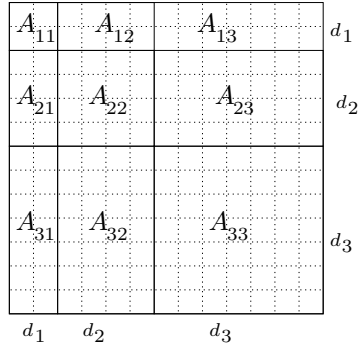}
\end{center}
\caption{The block decomposition of the matrix $A$ for $d=(2,4,7)$}\label{fig:blocks}
\end{figure}

Let $d=(d_1,\dots,d_t)$ be a dimension vector, $\frn$ the nilradical of the 
corresponding parabolic subalgebra. 
For $A\in\frn$ and 
$1\le i,j\le t$ 
we write $A_{ij}$ to describe the matrix formed by taking the entries of $A$ lying 
in the rectangle formed by rows $d_1+\dots + d_{i-1}+1$ up to 
$d_1+\dots + d_i$ and columns $d_1+\dots + d_{j-1}+1$ up to 
$d_1+\dots + d_j$ and with zeroes everywhere else. 
For $i\ge j$, this is just the zero matrix. 
Figure~\ref{fig:blocks} shows the blocks $A_{ij}$ for $d=(2,4,7)$. 

We set $A[i,j]$ to be 
the matrix formed by 
the $(A_{kl})_{i\le k\le j,i\le l\le j}$, 
i.e. by the rectangles right to and below of $A_{ii}$ and left to and 
above of $A_{jj}$. For instance, 
$A[i,i]$ is just $A_{ii}$ and 
$A[1,t]$ has the same entries as $A$. 
More generally, 
$A[ij]$ is a square matrix of size 
$(d_i+\dots+d_j)\times (d_i+\dots+d_j)$ with $A_{ii},\dots,A_{jj}$ on 
its diagonal. 

We are now ready to explain the rank conditions. For the rest of this 
section, we will always assume that a pair $(i,j)$ satisfies $1\le i<j\le t$.  
We write $X(d)$ for an element of $\cO(d)$.
For $k\ge 1$ define 
\begin{eqnarray*}
r_{ij}^k & := & \rk (X(d)[i,j]\,^k) \\
\kappa(i,j) & := & 1 + \#\{l\mid i<l<j,\ d_l\ge\min(d_i,d_j)\}\, .
\end{eqnarray*}
Observe that the 
numbers $r_{ij}^k$ are independent of the choice of an element of the 
Richardson orbit. 
With this, we can define two subsets of $\frn$ as our candidates 
for irreducible components of $Z$. 

\begin{definition}\label{def:Z-ij}
Let $d=(d_1,\dots,d_t)$ be a dimension vector and $\frn$ the nilradical 
of the parabolic subgroup $P$ of $\GL_n$. We set 
\begin{eqnarray*}
Z_{ij}^k & := & \{A\in\frn\mid \rk A[ij]^k<r_{ij}^k\} \\
Z_{ij} & := & Z_{ij}^{\kappa(i,j)}
\end{eqnarray*}
to be the elements $A$ of $\frn$ for which the rank of $k$th power of the matrix 
$A[ij]$ is defective, respectively the $A$ for which the rank of the 
$\kappa(i,j)$th power is defective. 
\end{definition}

To any dimension vector $d=(d_1,\dots,d_t)$ we associate subsets 
$\Gamma(d)$ and $\Lambda(d)$ of the set 
$\{(i,j)\mid 1\le i<j\le t\}$. 
In Section~\ref{s:rank-cond} we will show that the complement $Z$ 
of the open dense orbit is the union of the sets $Z_{ij}$ for $(i,j)\in \Lambda(d)$. 

\begin{eqnarray*}
\Gamma(d) & := & 
   \{(i,j)\mid d_l<\min(d_i,d_j)\ \mbox{or}\  d_l>\max(d_i,d_j)\ \forall\ i<l<j\}\, ,Ê\\
\Lambda(d) & := & 
  \left\{ (i,j)\in \Gamma(d)\mid d_i=d_j\right\} \cup \\ 
     & & \left\{ (i,j)\in \Gamma(d)\mid d_i\ne d_j \mbox{ and }\right. \\
     & & \left. 
     \begin{array}{lcl}
     \quad\quad & (i) & d_k\le\min(d_i,d_j)\ \mbox{or}\ d_k\ge\max(d_i,d_j)\ \forall\ k\\
     & (ii) & d_k\ne d_j\  \mbox{for $k<i$} \\
     & (iii) & d_k\ne d_i\  \mbox{for $k>j$}  
     \end{array} \right\},Ê\\
%
\end{eqnarray*}


Let us describe the latter in words: 
For $(i,j)$ to be in $\Lambda(d)$, we require that the $d_l$ with $i<l<j$ 
are smaller than the minimum of $d_i$ and $d_j$ or larger than the maximum 
of them. 
Furthermore, the $d_k$ have to be smaller or larger than 
the minimum $\min(d_i,d_j)$ resp. the maximum $\max(d_i,d_j)$ (for all $k$) 
and, if $d_i\ne d_j$, then $d_i$ is different from $d_{j+1},\dots,d_t$ and 
$d_j$ is different from $d_1,d_2,\dots,d_{i-1}$. 
In general, $\Gamma(d)$ is different from $\Lambda(d)$ as we illustrate now.

\begin{ex}\label{ex:Lambda}
\begin{itemize}
\item[(a)]
If $d=(1,3,4,2)$ then $\Gamma(d)=\{(1,2),(2,3),(3,4),(2,4),(1,4)\}$ 
and $\Lambda(d)=\{(2,3),(2,4),(1,4)\}$. 
\item[(b)]
For $d=(1,2,3,2)$, $\Gamma(d)=\{(1,2),(2,3),(3,4),(2,4)\}$, 
$\Lambda(d)=\{(1,2),(2,4)\}$. 
\item[(c)]
If $d=(d_1,\dots,d_t)$ is increasing or decreasing, then \\ 
$\Gamma(d)=\Lambda(d)$ 
$=\{(1,2),(2,3),\dots,(t-1,t)\}$. 
\item[(d)] The fourth example will be our running example throughout the paper: 
If $d=(7,5,2,3,5,1,2,6,5)$ then we have 
$\Gamma(d)=\{(i,i+1) \mid 1\le i\le 8\}$ \\
$\cup\ \{(1,8),(2,4),$ 
$(2,5),(3,6),(3,7),(4,6),(4,7),
(5,7),(5,8),(5,9),(7,9)\}$ 
and $\Lambda(d)=\{(1,8),(2,5),(3,7),(5,9)\}$. 
\end{itemize}
\end{ex}

We claim that the irreducible components of $Z=\frn\setminus \cO(d)$
are the $Z_{ij}$ with $(i,j)$ from the parameter set $\Lambda(d)$: 

\begin{theorem} (Theorem~\ref{thm:Z-ij})
Let $d=(d_1,\dots,d_t)$ be a composition of $n$, $\lambda=\lambda(d)$ 
the partition of the Richardson orbit corresponding to $d$. 
Then 
$$
Z=\bigcup_{(i,j)\in\Lambda(d)}  Z_{ij}
$$
is the decomposition of $Z$ into irreducible components. 
\end{theorem}

For the second description of the irreducible components we 
let $T(d)$ be the unique Young tableau obtained by filling the 
Young diagram of $\lambda$ with $d_1$ ones, $d_2$ twos, etc. 
(for details, we refer to Subsection~\ref{ss:young-tab}). 
Now for each pair $(i,j)$ we write $s(i,j)$ for the last row of 
$T(d)$ containing $i$ and $j$ and we let 
$T(i,j)$ be the tableau obtained from $T(d)$ by removing 
the box containing the number $j$ from row $s(i,j)$ and inserting 
it at the next possible position in order to obtain 
another tableau. The tableau $T(i,j)$ corresponds to an irreducible 
component of the intersection of $\frn$ with a nilpotent $\GL_n$-orbit 
as is explained in Section~\ref{s:tableaux} (Proposition~\ref{prop:comps}). 
We write $\frn(T(i,j))\subseteq\frn$ for 
the irreducible component in $\frn\cap C(\mu(i,j)$ of tableau $T(i,j)$. 
We claim that they correspond to irreducible components of $Z$ exactly 
for the $(i,j)\in\Lambda(d)$.

\begin{theorem} (Corollary~\ref{cor:tableaux})
Let $d=(d_1,\dots,d_t)$ be a dimension vector, $\lambda=\lambda(d)$ 
the partition of the Richardson orbit corresponding to $d$. 
Then 
$$
Z=\bigcup_{(i,j)\in\Lambda(d)} \frn(T(i,j))
$$
is the decomposition of $Z$ into irreducible components. 
\end{theorem}

As a consequence, we obtain that $Z$ has at most $t-1$ irreducible 
components (cf. Corollary~\ref{cor:components}) and we can describe their 
codimensions in $\frn$ (Corollary~\ref{cor:codim}). To be more precise, 
if $d$ is increasing or decreasing or if all the $d_i$ are different, 
then $Z$ has $t-1$ irreducible components. In particular, this applies 
to the Borel case where $d=(1,\dots,1)$. 
An example with $t=9$ and where we only have four irreducible components 
is our running example, see Example~\ref{ex:T(d)}. 

Note that the techniques we use are similar to the ones of~\cite{bah} 
where we describe the complement to the generic orbit in a representation 
space of a directed quiver of type A$_t$. However, the indexing sets 
are different and cannot be derived from each other. 

The paper is organised as follows: in Section \ref{s:rank-cond} we explain 
how to obtain the rank conditions. We first describe line diagrams associated to 
a composition $d$ of $n$. Line diagrams will be used to describe elements 
of the corresponding nilradical $\frn$. In Subsection~\ref{ss:Lambda-Gamma} 
we prove that the elements of $\Lambda(d)$ give the irreducible components. 
For this, we 
show that if $(i,j)$ does not belong to $\Gamma(d)$ then the variety $Z_{ij}$ is 
contained in a union of $Z_{k_sl_s}$ for a subset of elements $(k_s,l_s)$ 
of $\Gamma(d)$ 
(Lemma~\ref{lm:not-Gamma}). 
Next, if $(i,j)$ is in $\Gamma(d)\setminus\Lambda(d)$, then we can find 
$(k,l)\in\Lambda(d)$ such that $Z_{ij}$ is contained in $Z_{kl}$ 
(Corollary~\ref{cor:not-Lambda}). 
In Section~\ref{s:tableaux}, we recall Young diagrams and their fillings. 
Then we consider Young tableaux associated to a composition $d$ of $n$ and a 
nilpotency class $\mu\le\lambda(d)$. 
In a next step, we consider Young 
tableaux $T(i,j)$ associated to the elements of the parameter set $\Lambda(d)$. 
To each of these tableaux $T(i,j)$ we associate an irreducible variety 
$\frn(T(i,j))$: 
It is defined as the irreducible component in $\mathfrak n \cap C(\mu(i,j))$ corresponding 
to the tableau $T(i,j)$.
%
The $\frn(T(i,j))$ are known to 
be irreducible by work of the second author, \cite{hhab}. 
By showing that $\frn(T(i,j))$ is equal to $Z_{ij}$  
from Section~\ref{s:rank-cond} for elements $(i,j)$ of the parameter set $\Lambda(d)$ 
we can complete the description of the 
complement of the Richardson orbit in $\frn$ into irreducible 
components. 


%
\section{Components via rank conditions}\label{s:rank-cond}
%

%
\subsection{Line diagrams}
%
Let $d=(d_1,\dots,d_t)$ be a dimension vector for a parabolic subalgebra 
of $\frgl_n$, $\frn$ the corresponding nilradical. 
We recall a pictorial way to represent elements of $\frn$ and in particular, to 
obtain an element of the Richardson orbit
$\cO(d)$.
This can be found in~\cite[Section 2]{bhrr}  
and in~\cite[Section 3]{ba1}. 
We draw $t$ top-adjusted 
columns of $d_1$, $d_2$, $\dots, d_t$ vertices. 
The vertices are connected using edges between vertices of 
different columns. 
If two vertices lie on the same height and there is no third vertex between them 
on that height then we call the two vertices {\em neighbors}. 
The {\em complete line diagram} for $d$, $L_R(d)$, is the diagram with 
horizontal edges between all neighbored vertices 
(as the second and the third diagram of Example~\ref{ex:line-diagrams}). 
A {\em line diagram $L(d)$} for $d$ is a diagram with arbitrary edges between 
different columns (possibly with branching). 
A collection of connected edges is called a {\em chain of edges} (see the example 
below). If no branching occurs in a line diagram then a chain consisting of 
$l$ edges connects $l+1$ vertices. In that case we can define the length 
of a chain: 
The {\em length of a chain of edges} in a line diagram (without branching) 
is the number of edges 
the chain contains. A chain of length $0$ is a vertex that is not connected 
to any other vertex. 

In Example~\ref{ex:line-diagrams}, we show two complete and a branched 
line diagram for $d=(3,1,2,4)$ resp. for $d=(3,1,6,1,2,5,4)$. 

\begin{ex}\label{ex:line-diagrams}
a) A line diagram with branching and the complete line diagram $L_R(d)$ 
for $d=(3,1,2,4)$ are here. To the right of the latter we give the lengths 
of the chains in the diagram. 
\begin{center}
\includegraphics[scale=.5]{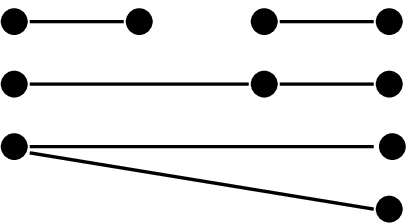} 
\hskip 20pt
\includegraphics[scale=.5]{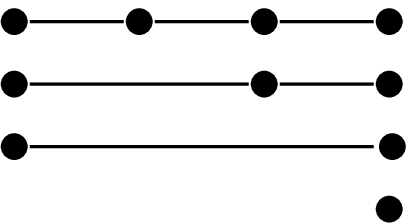} 
\end{center}
%


b) Now we consider our running example $d=(7,5,2,3,5,1,2,6,5)$. Its 
complete line diagram $L_R(d)$ is here, with the lengths of the chains to 
the right. 
\begin{center}
\includegraphics[scale=.5]{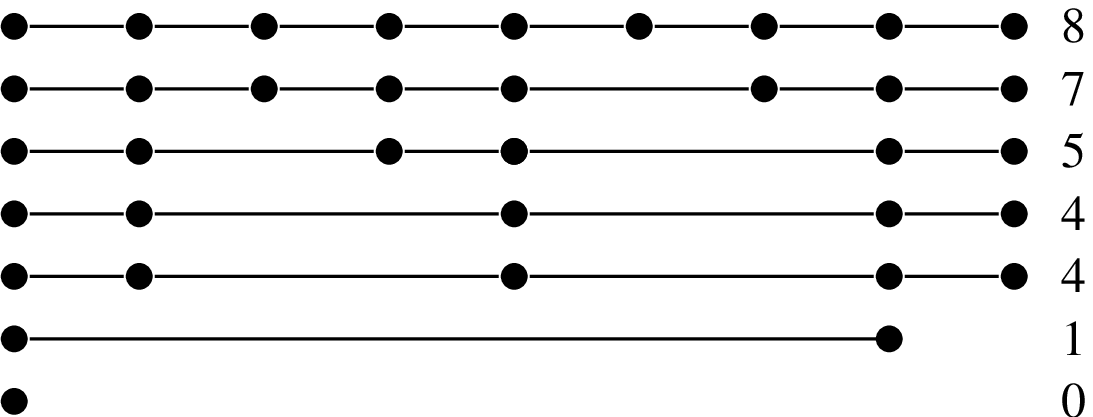}
\end{center}
\end{ex}

We will see in the next subsection 
that the line diagram $L_R(d)$ determines an element of the Richardson 
orbit of $\frn$. In general, 
line diagrams give rise to elements of 
the nilradical of nilpotency class smaller than $\lambda=\lambda(d)$ 
with respect 
to the Bruhat order. 
%


%


Any line diagram (complete or not) gives rise to an element $A$ of $\frn$: 

The sizes of the columns of a line diagram correspond to the sizes of the 
square blocks in the Levi factor of $\frp$. 
An edge between column $i$ and 
column $j$ (with $i<j$) of the diagram corresponds to a non-zero 
entry in the block $A_{ij}$ of the matrix $A$. 
A chain of two joint edges between three columns 
$i_0<i_1<i_2$ gives rise to a non-zero entry in block $A^2_{(i_0,i_2)}$ 
of the matrix $A^2$, etc. 
This can be made explicit, as we explain in the next subsection. 

%
\subsection{From line diagrams to the nilradical}
%

The elements of the nilradical $\frn$ for the dimension vector 
$d=(d_1,\dots, d_t)$ 
are nilpotent endomorphisms 
of $k^n$, for $n=\sum d_i$. 
In particular, if we write $e_1,\dots, e_n$ for a basis of $k^n$, then 
the elements of $\frn$ are sums $\sum_{i<j} a_{ij}E_{ij}$ for some 
$a_{ij}\in k$ where the elementary matrix $E_{ij}$ sends $e_j$ to $e_i$. 

We now describe a map associating an element of the nilradical to a given 
line diagram. We view the vertices of a line diagram $L(d)$ as 
labelled by the numbers $1,2,\dots, n$, 
starting at the top left vertex, with $1,2,\dots, d_1$ in the first column, 
$d_1+1,\dots,d_1+d_2$ in the second column, etc. 
Now if two vertices $i$ and $j$ (with $i<j$) are joint by an edge, we associate 
to this edge the matrix $E_{ij}$. 

We denote an edge between two vertices $i$ and $j$ 
($i<j\le n$) of the diagram by $e(i,j)$. Then we associate to an edge 
$e(i,j)$ of $L(d)$ the elementary matrix $E_{ij}\in\frn$. 
This can be extended to a map from the set of line diagrams for 
$d$ to the nilradical $\frn$ by linearity. 

For later use, we denote this map by $\Phi$: 
$$
\Phi: \{\mbox{line diagrams for $d$}\}\longrightarrow \frn, \ 
L(d)\mapsto \sum_{e(i,j)\in L(d)}E_{ij}\,.
$$

If $L(d)$ is a line diagram without branching, then the 
partition of the image under $\Phi$ of 
the line diagram $L(d)$ 
can be read off from it directly as follows: if $L(d)$ has $s$ chains of lengths 
$c_1, c_2, \dots, c_s$ (all $\ge 0$). 
Then $\sum_{j=1}^s (c_j+1)=\sum_{i=1}^t d_i=n$. 

\begin{remark}\label{re:line-part}
Let $L(d)$ be a line diagram without branching and let $c_1,\dots, c_s$ 
be the lengths of the chains of $L(d)$. 
Let $\mu=(\mu_1,\dots,\mu_s)$ be the partition obtained by ordering the 
numbers $c_j+1$ by size. Then $\mu$ is the partition of $\Phi(L(d))$. 
\end{remark}

In particular, $\Phi(L_R(d))$ is an element 
of the Richardson orbit $\cO(d)$
since the partition of $L_R(d)$ is just the dual of the dimension vector $d$ and this is 
equal to $\lambda(d)$ 
(cf. Section 3 in~\cite{ba1}). 
If $L(d)$ is any other line diagram for $d$ 
$L(d)$ (withouth branching), with lengths of chains $c_1,\dots,c_s$ and 
$\mu_i:=c_i+1$ 
then we always have $\sum_{j=1}^k\mu_j\le\sum_{j=1}^k\lambda_j(d)$ and 
so the partition of $\Phi(L(d))$ is smaller than or equal to the partition of 
$\Phi(L_R(d))$ under the Bruhat order. 
%

To summarize, we have the following: 
\begin{lemma}
Let $d$ be a dimension vector. Then, $\Phi(L(d))$ is an element of the 
nilradical $\frn$ of nilpotency class $\mu\le\lambda(d)$. 
In other words, $\Phi(L(d))$ lies in $\frn\cap C(\mu)$. 
\end{lemma}

\begin{ex}
Let $d=(3,1,2,4)$ as in Example~\ref{ex:line-diagrams} (a). The lengths of 
the chains of $L_R(d)$ are $3,2,1,0$, the Richardson orbit has partition 
$(4,3,2,1)$. 
We compute the matrix of the complete line diagram $L_R(d)$, 
and the powers of this matrix.  
Let $X(d):=\Phi(L_R(d))$. 
Then $X(d)$ and its powers are 
\begin{eqnarray*}
X(d)\  & = & E_{14}+E_{45}+E_{57}+E_{26}+E_{68}+E_{39}\\
X(d)^2 & = & E_{15}+E_{47}+E_{28} \\
X(d)^3 & = & E_{17} \\
X(d)^k & = & 0 \ \mbox{for $k>3$.}
\end{eqnarray*}

%
\end{ex}

Recall that we have defined the varieties $Z_{ij}^k$ by comparing 
the ranks of certain submatrices of elements in the nilradical $\frn$ to the 
corresponding rank $r_{ij}^k$ of a Richardson element, 
cf. Definition~\ref{def:Z-ij}. 
We thus need to be able 
to compute the rank of the submatrix $X(d)[ij]$ of an element 
$X(d)$ of the Richardson orbit $\cO(d)$ and of its powers.
For this, we can use the line diagram $L_R(d)$. 
Let $X(d)=\sum_{e(k,l)\in L_R(d)}E_{kl}$ be the Richardson 
element given by $L_R(d)$. 

To compute the rank $r_{1t}^k$ of $X(d)^k$, it is enough to count 
the chains of length $\ge k$ in the line diagram $L_R(d)$. 
Analogously, to find 
the rank $r_{ij}^k$ of the $k$th power of the submatrix 
$X(d)[ij]$, one has to count the chains of length 
$\ge k$ between the $i$th and $j$th column in $L_R(d)$: 

Let $1\le k<l\le n$ be such that the image $\Phi(e(k,l))$ of 
the edge $e(k,l)$ is in $X(d)[ij]$. 
That means we are considering edges $e(k,l)$ starting in some 
column $i_1\ge i$ and ending in some column $i_2\le j$. 
Thus, in computing $r_{ij}^k$, we really consider the $k$th 
power of the matrix which arises from columns $i,i+1,\dots,j$ 
of $L_R(d)$. We now introduce the notation to refer to the subdiagram 
consisting of these columns. 
We denote by $L_R(d)[ij]$ subdiagram of $L_R(d)$ of all vertices 
from the $i$th up to the $j$th column and of all edges starting strictly after the
$(i-1)$st column resp. ending strictly before the $(j+1)$st column. 
In other words, we remove columns $1,2,\dots,i-1$ and columns 
$j+1,\dots,t$ together with all edges incident with them. 

With this notation we have  
\begin{equation}\label{eq:s-r}
r_{ij}^k = \#\{\mbox{chains in $L_R(d)[ij]$ with at least $k$ edges}\}
\end{equation}
for $1\le i<j\le t$, $k\ge 1$.

Similarly, if $L(d)$ is a 
line diagram for $d$, we write $L(d)[ij]$ 
to denote the subdiagram of $L(d)$ of rows $i$ to $j$.

\begin{ex}
The subdiagram $L_R(d)[47]$ for $d=(7,5,2,3,5,1,2,6,5)$ of the diagram 
$L_R(d)$ from (b) of Example~\ref{ex:line-diagrams} is 
shown here (dotted lines and empty circles are thought to be 
removed): 
\begin{center}
\includegraphics[scale=.5]{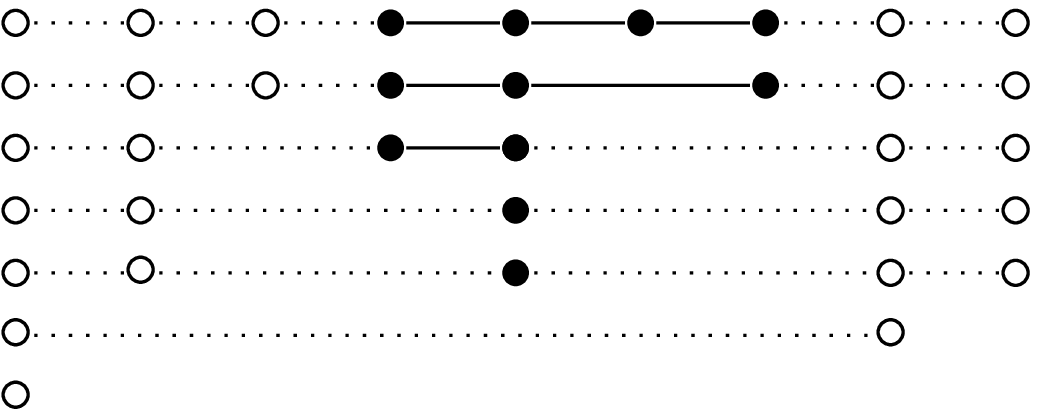}
\end{center}
\end{ex}

%
\subsection{The varieties $Z_{ij}$}\label{ss:Lambda-Gamma}
%

As explained earlier, we want to show that 
the irreducible components of $Z$ are indexed by 
the parameter set $\Lambda(d)$. 
With this in mind, we now 
discuss the properties of the varieties $Z_{ij}^k$. 
We will prove that for $l\ne\kappa(i,j)$, $Z_{ij}^l$ is either empty or contained in 
$Z_{ij}$ or in the union $Z_{ij_0}\cup Z_{i_0j}$ for some $i_0\le j_0$. 
Later in this section we will see that not all $(i,j)$ with $1\le i<j\le t$ are needed 
to describe the complement $Z$. 

The following notations will be useful: 
\begin{eqnarray*}
d_{<}[ij] & := & |\{l\mid i<l<j,\ d_l<\min(d_i,d_j)\}| \\ 
d_{\ge}[ij] & := & |\{l\mid i<l<j,\ d_l\ge\min(d_i,d_j)\}| \, . 
\end{eqnarray*} 
%
If $d=(7,5,2,3,5,1,2,6,5)$, then $d_<[25]=2$, $d_<[26]=\emptyset$ 
and $d_{\ge}[26]=3$. 

\begin{figure}[h]
\begin{center}
\includegraphics[scale=.5]{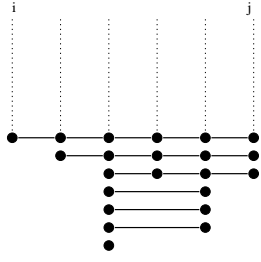}
\end{center}
\caption{Our running example has 
$d_{<}[26]=\emptyset$, so $\kappa(2,6)=4$}\label{fig:no-min}
\end{figure}

\begin{remark}\label{re:kappa}
Observe that 
\begin{eqnarray*}
\kappa(i,j) & = & 1+ \# d_{\ge}[ij] \\
  & = & j-i-\# d_{<}[ij]\,. 
\end{eqnarray*}
In particular, 
$\kappa(i,j)=j-i$ if and only if $d_{<}[ij]=\emptyset$. Figure~\ref{fig:no-min}Ê
illustrates this. 
\end{remark}

\begin{lemma}\label{lm:emptyset}
Let $d=(d_1,\dots,d_t)$ be a dimension vector and $1\le i<j\le t$. 
Then for $k>0$ we have 
$$
Z_{ij}^k=\emptyset \mbox{ if and only if }
k> j-i \,.
$$ 
\end{lemma}

\begin{proof}
One has $r_{ij}^k=\rk X(d)[ij]^k>0$ exactly for $k\le j-i$ and 
$0\in Z_{ij}^k$ if and only if $r_{ij}^k>0$. 
\end{proof}

It remains to consider the cases where $l$ is smaller than $\kappa(i,j)$ 
or when $l$ lies between $\kappa(i,j)$ and $j-i$. 
This is covered by the next two statements. 

\begin{lemma}\label{lm:l<kappa}
For $1\le l< \kappa(i,j)$ the following holds: 
$$
Z_{ij}^l\subsetneq Z_{ij}\,.
$$
\end{lemma}

\begin{proof}

We may assume $d_i\le d_j$. 
For any $B\in\frn$ the rank of $B[ij]^l$ is independent of the order of 
$d_i,d_{i+1},\dots,d_j$: incomputing the rank, we need to know the 
number of (independent) chains of length $l$ in 
the line diagram of $b[ij]$. Hence we may reorder $d_i,\dots,d_j$ 
to obtain $d_{s_1},\dots, d_{s_{j-i+1}}$ with 
$d_{s_k}\le d_{s_{k+1}}$ 
for $k=1,\dots,j-i$.
One computes $r_{ij}^l=\rk X(d)[ij]^l$ as the 
sum $\sum_{k=0}^{j-i-l}d_{i+k}$. 

Let $A$ belong to $Z_{ij}^l$ for some $l<\kappa(i,j)$. 
Thus $\rk A[ij]^l<r_{ij}^l=\rk X(d)[ij]^l$. 
But then also the rank of $A[ij]^k$ is smaller than $r_{ij}^k$ for 
$k=l+1,\dots,\kappa(i,j)$. 
In particular, $A\in Z_{ij}$. 
The inequality is clear. 
\end{proof}

Let $A$ belong to $Z_{ij}^l$ for some $l<\kappa(i,j)$. 
Thus $\rk A[ij]^l<r_{ij}^l$. 
But then also the rank of $A[ij]^k$ is smaller than $r_{ij}^k$ for 
$k=l+1,\dots,\kappa(i,j)$. 
In particular, $A\in Z_{ij}$. 

\begin{lemma}\label{lm:l>kappa}
For $\kappa(i,j)<l\le j-i$ 
the following holds: there exist $i_0\le j_0\in d_{<}[ij]$, 
$d_{i_0}$, $d_{j_0}< \min(d_i,d_j)$ maximal, such that 
$$
Z_{ij}^l\subseteq Z_{ij_0} \cup Z_{i_0j}\,.
$$
\end{lemma}

\begin{proof}
We first observe that for elements of the Richardson orbit, the rank 
%
%
%
$r_{ij}^l$ is 
$$
r_{ij}^l=\sum_{i_0=i}^{j-l}\max_{\begin{tiny}
                       \begin{array}{c}
                                                 i_0<\dots<i_{l}\le j
                       \end{array}
                       \end{tiny}}
                       \min\{d_{i_0},\dots,d_{i_l}\}
$$

(1) 
Let us first consider the case where $d_{<}[ij]$ only has one element, say 
$d_{<}[ij]=\{i_0\}$, see Figure~\ref{fig:one-min}). 
Then $\kappa(i,j)=j-i-1$ and so $l=j-i$. 

%

For $A\in\frn$ to be an element of $Z_{ij}^{l}$, 
the rank of $A[ij]^l$ is smaller than $r_{ij}^l$. Since $d_{i_0}$ is 
minimal among all $d_i,\dots,d_j$, this implies $\rk A[ii_0]^l<r_{ij}^l$ or 
$\rk A[i_0j]^l<r_{ij}^l$ and we are done. 

\begin{figure}[h]
\begin{center}
\includegraphics[scale=.5]{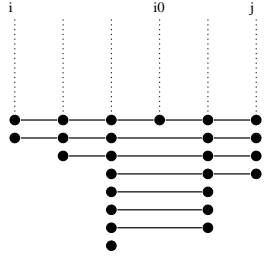}
\end{center}
\caption{The case $|d_<[ij]|=1$: in the running example, we have 
$d_{<}[47]=\{6\}$ and $\kappa(i,j)=2$}\label{fig:one-min}
\end{figure}

(2)
The case where $d_{<}[ij]$ has at least two elements only 
needs a slight modification of the argument. 
Take $i_0$, $j_0$ from $d_{<}[ij]$ with $d_{i_0}$, $d_{j_0}$ maximal with 
$i_0$ being the smallest among these indices, $j_0$ the largest one 
(we do not distinguish between the two possibilities 
$d_{i_0}=d_{j_0}$ and $d_{i_0}\ne d_{j_0}$), see Figure~\ref{fig:more-min}. 
With a similar reasoning as in part (1) of the proof, $A$ then lies in 
$Z_{i,j_0}$ or in $Z_{i_0,j}$. 

\begin{figure}[h]
\begin{center}
\includegraphics[scale=.5]{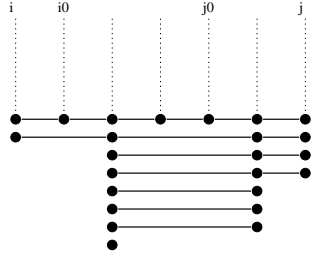}
\end{center}
\caption{The case $i_0\ne j_0\in d_{<}[ij]$: our running example has 
$d_<[48]=\{6,7\}$ and so $\kappa(4,8)=2$}\label{fig:more-min}
\end{figure}

\end{proof}


\begin{lemma}\label{lm:Z-union}
The complement $Z$ decomposes as follows:
$$
Z=\cup_{1\le i<j\le t} Z_{ij}= \cup_{ij} \cup_{k\ge 1} Z_{ij}^k\, .
$$
\end{lemma}

\begin{proof}
The inclusion $\subseteq$ of the second equality is clear. 
To obtain the inclusion $\supseteq$, one uses 
Lemmata~\ref{lm:emptyset},~\ref{lm:l<kappa} and~\ref{lm:l>kappa}. 
Consider the first equality: by definition, $A\in Z$ if and only 
if $A\notin\cO(d)$. The latter is the case if and only if there exist
$1\le i<j\le t$, $k\le j-i$, such that $A\in Z_{ij}^k$: to see this, one uses 
the formula for the dimension of the stabilizer of $A\in \frgl_n$, see \cite{kp}. 
This formula uses the dimensions of the kernels of the maps $A^k$, $k\ge 1$. 
The stabilizer of $A$ has dimension $0$ if and only if $A$ is an element of $\cO(d)$. 
\end{proof}

It now remains to see that the $(i,j)\in\Lambda(d)$ are enough to 
describe the irreducible components of $Z$. In a first step 
(Lemma~\ref{lm:not-Gamma}), we start with 
$(i,j)\notin \Gamma(d)$ and show that in that case 
$Z_{ij}$ is contained in a union of $Z_{kl}$'s such that the corresponding $(k,l)$ 
all lie in $\Gamma(d)$. 

Then we consider an element $(i,j)$ of 
$\Gamma(d)\setminus\Lambda(d)$ and 
show that we can find $(k,l)$ $\in\Lambda(d)$ 
with $Z_{ij}\subseteq Z_{kl}$ (Lemma~\ref{lm:not-Lambda} 
and Corollary~\ref{cor:not-Lambda}). As always, we assume that 
$1\le i<j\le t$ and $1\le k<l\le t$. 

\begin{lemma}\label{lm:not-Gamma}
Assume that $(i,j)$ does not belong to $\Gamma(d)$. 
Then there exists $\Gamma'(d)\subseteq\Gamma(d)$ such that 
$$
Z_{ij}\subseteq \bigcup_{(k,l)\in\Gamma'(d)} Z_{kl}\,.
$$
\end{lemma}

\begin{proof}
It is enough to show that we can find an $l$, $i<l<j$, with 
$\min(d_i,d_j)\le d_l\le \max(d_i,d_j)$, such that 
$$
Z_{ij}\subseteq Z_{il}\cup Z_{lj}\, .
$$
By iterating this, we will eventually end up with 
a subset $\Gamma'(d)\subset\Gamma(d)$ as in the statement of the lemma. 
%
%

So choose an $l$, $1<l<t$, with $\min(d_i,d_j)\le d_l\le \max(d_i,d_j)$ 
(such an $l$ exists since $(i,j)\notin\Gamma(d)$). Take $A\in Z_{ij}$ arbitrary. 
By assumption, $A[ij]^{\kappa(i,j)}$ is defective, i.e. 
$\rk A[ij]^{\kappa(i,j)}<r_{ij}^{\kappa(i,j)}$. Since $d_l\ge d_i,d_j$, the defectiveness 
is inherited from $A[il]$ or from $A[lj]$ and $A\in Z_{il}$ or $A\in Z_{lj}$ accordingly. 
\end{proof}

\begin{figure}[h]
\begin{center}
\includegraphics[scale=.5]{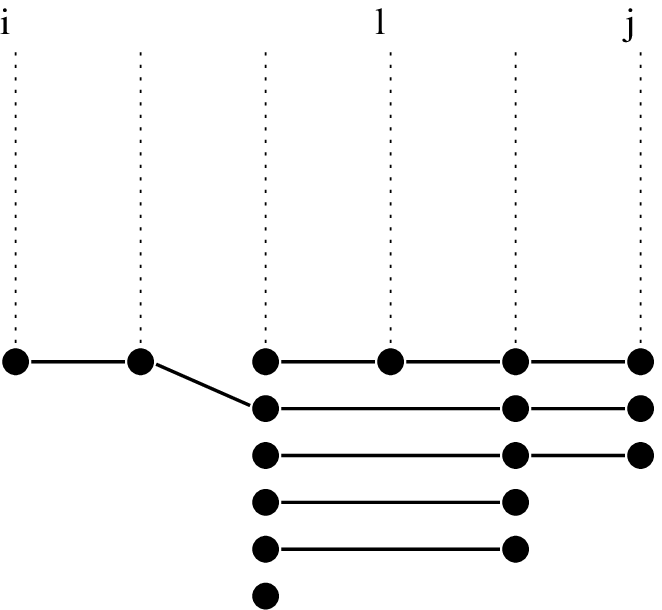} \hspace{1cm}
\includegraphics[scale=.5]{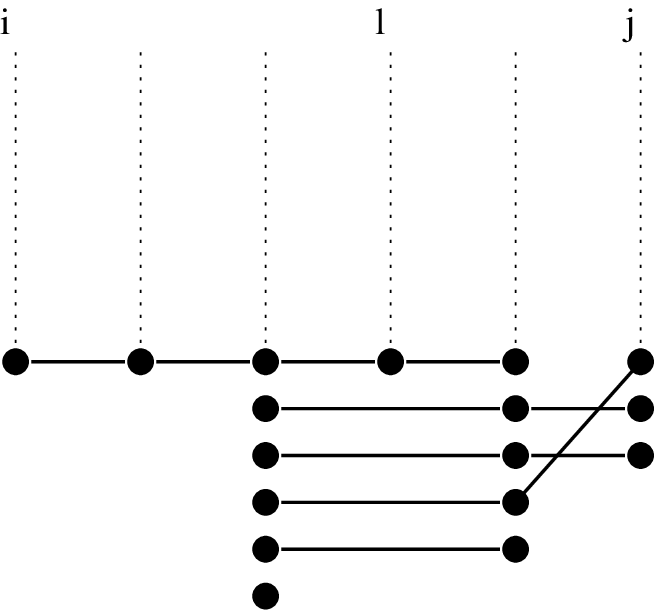}
\end{center}
\caption{Examples for $A\in Z_{il}$ resp. for $A\in Z_{lj}$ for $d=(7,5,2,3,5,1,2,6,5)$, 
with $i=3$, $j=8$ and $l=6$}\label{fig:not-Gamma}
\end{figure}

%
Let us remark that when removing an edge of a chain of $L_R(d)$ 
in the proof above, we ensured that the matrix $A$ has a zero entry 
at the corresponding position. 
In general, the diagram of a matrix in $Z_{il}$ resp. in $Z_{lj}$ 
has more non-zero entries than the ones obtained after removing 
one edge from $L_R(d)$: this is illustrated by the dashed lines in 
Figure~\ref{fig:not-Gamma}. 

The following lemma states that for any $(i,j)$ from $\Gamma(d)\setminus\Lambda(d)$ 
there exists $(k,l)$ from $\Lambda(d)$ with 
$k\le i<j\le l$ such that $Z_{ij}\subseteq Z_{kl}$. 

\begin{lemma} \label{lm:not-Lambda}
Assume that $(i,j)\in\Gamma(d)\setminus \Lambda(d)$. 
Then one of the following holds: 
$$
\begin{array}{ll}
  & \mbox{there exists $k>j$ with $Z_{ij}\subseteq Z_{ik}$} \\
\mbox{or} & \mbox{there exists $l<i$ with $Z_{ij}\subseteq Z_{lj}$.}
\end{array}
$$
\end{lemma}


\begin{proof}
First observe that $d_i\ne d_j$ since $(i,j)$ belongs to $\Lambda(d)$ otherwise. 
Without loss of generality, we assume $d_i<d_j$. 
We have three cases to consider: 
\begin{itemize}
\item[(i)]
There is $k_1\in\{1,\dots,i-1\}\cup\{j+1,\dots,t\}$ with $d_i<d_{k_1}<d_j$. 
\item[(ii)]
There exists $k_2<i$ with $d_{k_2}=d_j$. 
\item[(iii)]
There exists $k_3>j$ with $d_{k_2}=d_i$. 
\end{itemize}
\begin{figure}[h]
\begin{center}
  \psfragscanon
   \psfrag{1}{$m_1$}
   \psfrag{2}{$m_2$}
   \psfrag{4}{$m_4$}
   \psfrag{5}{$m_5$}
   \psfrag{9}{$m_9$}
\includegraphics[scale=.6]{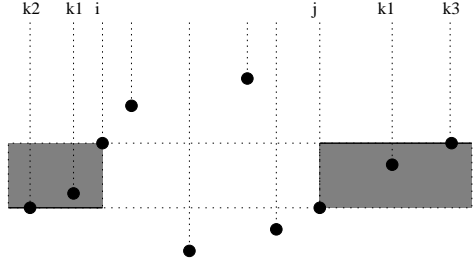}
\caption{For the running example, $(7,8)$ is in $\Gamma(d)\setminus\Lambda(d)$, 
as for all $i$, $d_{m_i}$ violates the assumptions on $\Lambda(d)$}\label{fig:not-Lambda}
\end{center}
\end{figure}
The three cases are illustrated in Figure~\ref{fig:not-Lambda}: if $(i,j)\in\Gamma(d)$ 
but not in $\Lambda(d)$ then one of the following has to occur: there has to be a 
$k$ with $d_k$ inside the shaded area or with $d_k$ lying on the same height as 
$d_j$ (if $k<i$) resp. on the same height as $d_i$ (if $k>j$). 


\noindent
Case (i) with $k_1>j$: Among the $k_1>j$ with $d_i<d_{k_1}<d_j$ choose one with 
$d_{k_1}-d_i$ minimal, and $k_1$ minimal (i.e. as close to $j$ as possible). 
Note that we have $\kappa(i,j)\le \kappa(i,k_1)$. 
Now $A\in Z_{ij}$ means that $A[ij]^{\kappa(i,j)}$ is defective. 
Since $d_i<d_{k_1}$, this defectiveness has to be inherited from 
$A[i,k_1]$, i.e. $\rk A[i,k_1]^{\kappa(i,k_1)}<r_{i,k_1}^{\kappa(i,k_1)}$ and 
so, $Z_{ij}\subseteq Z_{i,k_1}$.  


\noindent
Case (i) with $k_1<i$: here, we choose $k_1$ accordingly to be such that 
$d_j-d_{k_1}$ is minimal and $k_1<i$ maximal among those (i.e. as close 
to $i$ as possible). One checks that $\kappa(i,j)\le\kappa(k_1,j)$. Similarly 
as before, one gets $Z_{ij}\subseteq Z_{k_1,j}$. 


\noindent
Case (ii) : Among the $k_2<i$ with $d_{k_2}=d_j$, choose the maximal one 
(i.e. the one closest to $i$). We have $\kappa(i,j)\le\kappa(k_2,j)$ and 
we get $Z_{ij}\subseteq Z_{k_2,j}$. 
Case (iii) is completely analogous to case (ii). 
\end{proof}

Observe that $(k_2,j)$ and $(i,k_3)$ from cases (ii) and (iii) above 
are elements of $\Lambda(d)$.

\begin{cor}\label{cor:not-Lambda}
For any $(i,j)\in\Gamma(d)\setminus\Lambda(d)$ there exists 
$(k,l)\in\Lambda(d)$ such that 
$$
Z_{ij}\subseteq Z_{kl}\, .
$$
\end{cor}

\begin{proof}
Without loss of generality, we can assume $d_i<d_j$. 
By the observation after the proof of Lemma~\ref{lm:not-Lambda}, we are 
done if there exists $k'<i$ with $d_{k'}=d_j$ or $k''>j$ 
with $d_{k''}=d_i$. Using similar arguments, one sees that if there exist 
$k'<i$ and $k''>j$ with $d_i <d_{k'}=d_{k''}< d_j$ 
then $(k',k'')\in\Lambda(d)$ and $Z_{ij}\subseteq Z_{k',k''}$. 
Thus, assume that there exists 
$k\in\{1,\dots,i-1\}\cup \{j+1,\dots,t\}$ with $d_i <d_k< d_j$ 
and such that there is no $k'<i$ with $d_{k'}=d_j$ and no $k''>j$ 
with $d_{k''}=d_i$. \\
If $k>j$, we choose $k$ such that $d_k-d_i$ is 
minimal and take the minimal $k>j$ among these 
(i.e. $k$ is as close to $j$ as possible). 
There are two possibilities: \\
Either we have $d_{k'}>d_k$ for all $k'<i$. Then, 
$(k',k)\in\Lambda(d)$ and one checks that $Z_{ij}\subseteq Z_{k',k}$. \\
Or there exists is $k'<i$ with $d_i<d_{k'}<d_k$. In that case, among the $k'<i$ 
with this property, we choose one with $d_k-d_{k'}$ minimal and such that 
$k'<i$ is maximal (i.e. $k'$ is as close to $i$ as possible). Again, we get 
$(k',k)\in\Lambda(d)$ and $Z_{ij}\subseteq Z_{k',k}$. \\
The case $k<i$ is analogous.  
\end{proof}



%
\section{Components via tableaux}\label{s:tableaux}
%

Let $d=(d_1,\dots,d_t)$ be a composition of $n$ and $\cO(d)$ be
the corresponding Richardson orbit in $\frn$, let $\lambda=\lambda(d)$ 
be the partition of the Richardson orbit. 
The second description of the irreducible components of $Z=\frn\setminus\cO(d)$
uses partitions $\mu_{ij}$, for $(i,j)\in\Lambda(d)$ and tableaux corresponding 
to them. 
Observe that $\lambda_1=t$, that $\lambda_2$ is the number of $d_i\ge 2$ 
appearing in $d$, 
$\lambda_3=\#\{d_i\mid d_i\ge 3\}$, and so on.

Let us introduce the necessary notation. 
If $\lambda=\lambda_1\ge\lambda_2\ge\dots\ge \lambda_s\ge 1$ is a partition 
of $n$ we will also use $\lambda$ to 
denote the Young diagram of shape $\lambda$. It has 
$s$ rows, with $\lambda_1$ boxes in the top row, $\lambda_2$ boxes in the 
second row, etc., up to $\lambda_s$ boxes in the last row. That means 
that we view Young diagrams as a number of right adjusted rows of boxes, 
attached to the top left corner, and decreasing in length from top to 
bottom. 
A standard reference for this is the book~\cite{fu} 
by Fulton. 
%

%
\subsection{The Young tableaux $\cT(\mu,d)$}\label{ss:young-tab}
%

Let $\mu\le\lambda(d)$ be a partition of $n$ (unless mentioned otherwise, we 
will always deal with partitions of $n$).  

\begin{definition}
We define a {\em Young tableau} 
of shape $\mu$ and of dimension vector $d$ 
to be a  filling 
of the Young diagram of $\mu$ with $d_1$ ones, $d_2$ twos, etc.
We write $\cT(\mu,d)$ for the set of all Young tableaux 
of shape $\mu$ and for $d$. 
\end{definition}

Recall that the rules for fillings of a Young diagram are that the numbers in 
a row strictly increase from left to right and that the numbers in a column 
weakly increase from top to bottom. In general, there might be several 
Young tableaux of a given shape for a given $d$. 
There is exactly one Young tableau of shape 
$\lambda=\lambda(d)$ and for 
$d$, so $\cT(\lambda(d),d)$ only has one element. 
To abbreviate, we will just call it $T(d)$. The entries of the boxes of its first row 
are $1,2,\dots,t$.

\begin{ex}
The partition of the composition $d=(7,5,2,3,5,1,2,6,5)$ of 
36 is $\lambda(d)=(9,8,6,5,5,2,1)$. The partition $\mu=(9,8,6,5,4,3,1)$ 
is smaller than $\lambda(d)$ and $\cT(\mu,d)$ consists of one element 
$T(\mu,d)$. We include $T(d)$ and $T(\mu,d)$ here.  

$$
\psfragscanon
   \psfrag{1}{$_{1}$}
   \psfrag{3}{$_{3}$}
   \psfrag{2}{$_{2}$}
   \psfrag{4}{$_{4}$}
   \psfrag{5}{$_{5}$}
   \psfrag{6}{$_{6}$}
   \psfrag{7}{$_{7}$}
   \psfrag{8}{$_{8}$}
   \psfrag{9}{$_{9}$}
T(d)\ 
\includegraphics[scale=.5]{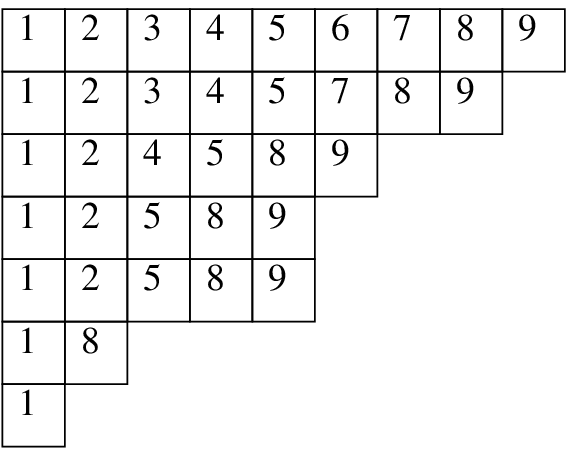} \hspace{1cm}
T(\mu,d)\ 
\includegraphics[scale=.5]{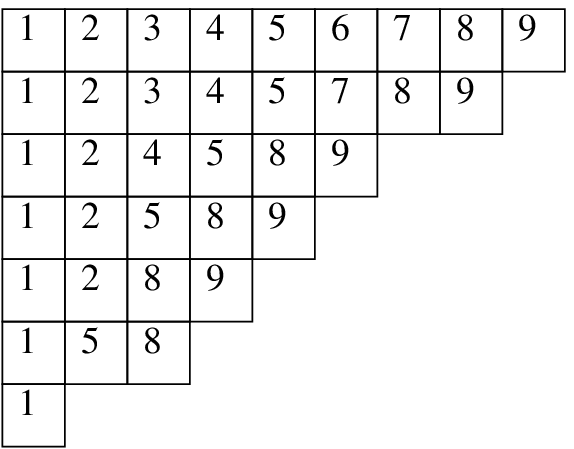} 
$$
\end{ex}

In order to understand the irreducible components of 
the complement $Z=\frn\setminus \cO(d)$, we have to consider the
intersections $\frn\cap C(\mu)$ for $\mu<\lambda(d)$. 
Each irreducible component of $Z$ corresponds to 
an irreducible component in such an intersection. 
Here, we can use a result of the second author (cf. Section 4.2 of~\cite{hhab}). 
First, one observes that the irreducible components of $\frn\cap C(\mu)$ 
are given by sequences $\mu^1,\dots,\mu^t$ where $\mu^i$ is a partition of 
$\sum_j^i d_j$ where $\mu^t=\mu$ and such that 
$0\le\mu_j^{i+1}-\mu_j^i\le 1$ (for all $j$, for $1\le i<t$). And the latter 
correspond to tableaux of shape $\mu$ with $d_i$ entries $i$, i.e. the elements 
of $\cT(\mu,d)$ in our notation.

\begin{prop}\label{prop:comps}
Let $\mu\le \lambda(d)$ be a partition of $n$. Then the 
irreducible components of $\frn\cap C(\mu)$ are in natural 
bijection with 
with the tableaux in $\cT(\mu,d)$. 
\end{prop}

\begin{proof}
This is Satz 4.2.8 in~\cite{hhab}.
\end{proof}

\begin{ex} \label{ex:O(d)-T(d)}
Let $d=(d_1,\dots,d_t)$ be a dimension vector and $\lambda=\lambda(d)$. 
We know that $\frn\cap C(\lambda)=\cO(d)$ is the Richardson orbit.
On the other hand, $\cT(\lambda,d)=T(d)$ has exactly one tableau. 
We now explain how to relate the complete line diagram $L_R(d)$ to 
the tableau $T(d)$. 
The lengths of the chains in $L_R(d)$ are the entries of the partition of $\lambda$ and hence 
give the shape of $T(d)$. 
The filling of $T(d)$ can now be obtained from $L_R(d)$ by labelling each vertex of the $i$-th column 
in $L_R(d)$ by an $i$. These numbers are then copied row by row, from left to right 
into the Young diagram of shape $\lambda$ to get $T(d)$. 


$$
\psfragscanon
   \psfrag{1}{$_{1}$}
   \psfrag{3}{$_{3}$}
   \psfrag{2}{$_{2}$}
   \psfrag{4}{$_{4}$}
   \psfrag{label columns}{\tiny{\mbox{label columns}}}
   \psfrag{copy into T(d)}{\tiny{\mbox{copy into $T(d)$}}}
\includegraphics[scale=.6]{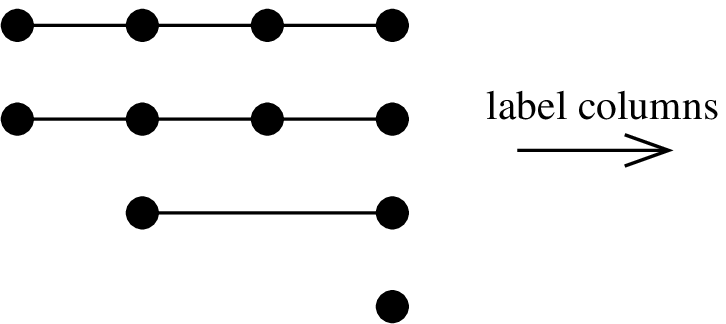}\quad\quad
\includegraphics[scale=.6]{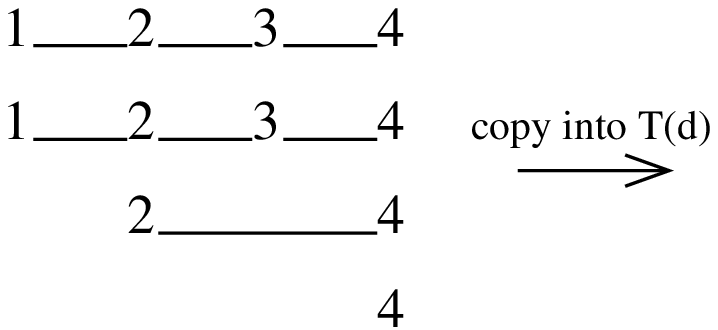}\quad\quad
\includegraphics[scale=.6]{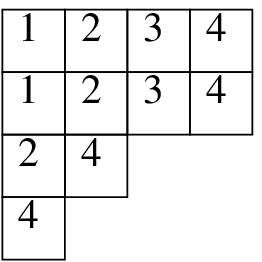}
$$
\end{ex}

From this connection between the line diagram $L_R(d)$ and $T(d)$ one 
deduces the following useful observation. 
Every pair $(i,j)$ with $1\le i<j\le t$ determines a unique row of $T(d)$ namely
the last row of $T(d)$ containing
$i$ and $j$.
Such a row always exists as the first row just consists of the boxes
with numbers $1,2,3,\dots,t$.
We denote this row by $s(i,j)$.

\begin{lemma}\label{lm:kappa-T(d)}
The number of boxes between $i$ and $j$ 
in row $s(i,j)$ of $T(d)$ is equal to $\kappa(i,j)-1$.
\end{lemma}

Proposition~\ref{prop:comps} describes the irreducible components of 
the intersections $\frn\cap C(\mu)$ for $\mu\le \lambda$: 
They are given by the Young tableaux in $\cT(\mu,d)$, i.e. 
by all possible fillings of the diagram $\mu$ by the numbers given by $d$.


Clearly, not all irreducible components of the different
intersections $\frn\cap C(\mu)$ give rise to an
irreducible component of $Z$.
If $\mu_2\le \mu_1$
and $T_i\in \cT(\mu_i,d)$ are tableaux such that
$T_2$ can be obtained from $T_1$ by moving down boxes successively,
then the irreducible component corresponding to $T_2$ is already
contained in the irreducible component corresponding to $T_1$
and thus does not give rise to a new irreducible component of the complement
$Z$ of the Richardson orbit.
This is in particular the case, if $T_1$ is obtained from the tableau $T(d)$
of the Richardson orbit by moving down a single box
and $T_2$ is a degeneration of $T_1$ (obtained by moving down
boxes from $T_1$).
Thus, the only candidates for irreducible components are the ones
given by tableaux which can be obtained
from $T(d)$ by moving down a single box 
to the closest possible row. 
We call such a degeneration a {\em minimal movement}.

%
\subsection{The Young tableaux $T(i,j)$}
%

To describe minimal movements, we now define certain tableaux $T(i,j)$.

\begin{definition}
The tableau $T(i,j)$ is the tableau obtained from $T(d)$ by removing 
the box containing the number $j$ from row $s(i,j)$ and inserting it 
in the nearest row in order to obtain another tableau. 
In other words: Among the possible rows where this box could be inserted, 
we choose the one that is closest to row $s(i,j)$. 
We denote the partition of the resulting tableau $T(i,j)$ by $\mu(i,j)$. 
\end{definition}

\begin{definition}
For a tableau $T(i,j)$ we define $\frn(T(i,j))\subseteq \frn$ to be the irreducible 
component of $\frn\cap C(\mu(i,j))$ whose tableau is $T(i,j)$. 
%
\end{definition}
We claim that $\frn(T(i,j))$ gives rise to an irreducible component of the complement 
$Z$ exactly when 
$(i,j)$ belongs to the parameter set $\Lambda(d)$. 

For completeness, we recall the definition of a the tableau $T$ 
for a an irreducible component in
$C(\mu) \cap \mathfrak n$. Consider a maximal flag $V_0 \subset V_1 \subset
\ldots \subset V_t$ of vector spaces that is stabilized by $P(d)$. Take any
matrix $A$ in the open subset of an irreducible component of $C(\mu) \cap
\mathfrak n$ where $A$ restricted to $V_i$ has constant Jordan type. Then the
Young diagram of $A|_{V_{i}}$ is the partition obtained from $T$ by deleting
all boxes with entries $i+1,\ldots,t$. So the subdiagramm consisting of all
boxes with entries at most $i$ measures the generic Jordan type of $A$
restricted to the subspace $V_i$. In particular, the equation defining the
component corresponding to $T(i,j)$ can involve only equations in the entries
of $A[1,j]$. Even stronger, we will see in Lemma~\ref{lm:T-gleich-Z} 
that the equations
involve only entries in $A[i,j]$ for $(i,j)$ in $\Gamma(d)$.

\medskip

To prepare for Lemma~\ref{lm:T-gleich-Z} we observe 
that for $(1,t)\in\Gamma(d)$ the component $\frn(T(1,t))$ 
coincides with $C(\mu(1,t)) \cap \mathfrak n$ since there is only one tableau
for the partition $\mu(1,t)$ with dimension vector $d$. Consequently, this
component is defined by the equation 
$\rk A[1,t]^{\kappa(1,t)} <r_{1,t}^{\kappa(1,t)}$ defining $C(\mu(1,t))\cap\frn$ inside $\frn$.

\medskip


By definition, the tableau $T(i,j)$ is 
obtained from $T(d)$ through a minimal movement. 
Its partition $\mu(i,j)$ is clearly smaller than $\lambda=\lambda(d)$ as 
the lengths of the rows of a tableau are the parts of 
the corresponding partition. In particular, these lengths 
form a decreasing sequence of positive numbers. 
Thus, moving down a box from a row of length $k$ to a lower row 
of length at most $k-2$ 
results in a partition which is smaller than the 
original partition. 
Note, however, that different elements $(i,j)$ and $(k,l)$ can lead to the 
same partition $\mu(i,j)=\mu(k,l)$, e.g. 
$\mu(2,5)=\mu(5,9)$ in Example~\ref{ex:T(d)} below. 

\begin{ex}\label{ex:T(d)}
Let $d=(7,5,2,3,5,1,2,6,5)$ be a dimension vector, $n=36$. 
To illustrate the construction of $T(i,j)$ we compute these tableaux for all 
$(i,j)$ $\in\Lambda(d)=\{(1,8),(2,5),(3,7),(5,9)\}$. They are presented in 
Figure~\ref{fig:T(d)}. In the picture showing the line diagram $L_R(d)$ we have 
indicated the connections between the columns $i$ and $j$ for all 
$(i,j)\in\Lambda(d)$ by shaded areas. 
\begin{figure}[ht]
\begin{center}
  \psfragscanon
   \psfrag{1}{$_{1}$}
   \psfrag{2}{$_{2}$}
   \psfrag{3}{$_{3}$}
   \psfrag{4}{$_{4}$}
   \psfrag{5}{$_{5}$}
   \psfrag{6}{$_{6}$}
   \psfrag{7}{$_{7}$}
   \psfrag{8}{$_{8}$}
   \psfrag{9}{$_{9}$}
$T(d)$
\includegraphics[scale=.5]{T_d.eps} \\[.4cm]
$T(1,8)$
\includegraphics[scale=.5]{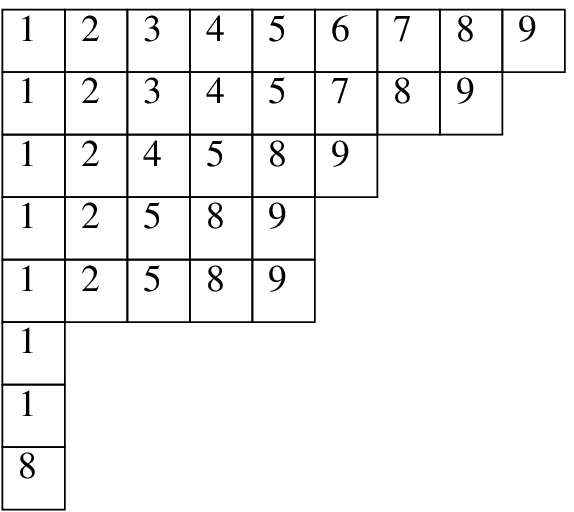} \vspace{.4cm}
$T(2,5)$
\includegraphics[scale=.5]{T_25.eps} \\
$T(3,7)$
\includegraphics[scale=.5]{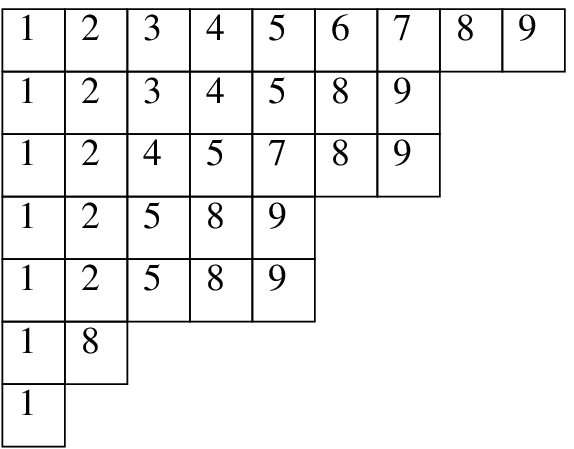} \vspace{.4cm}
$T(5,9)$
\includegraphics[scale=.5]{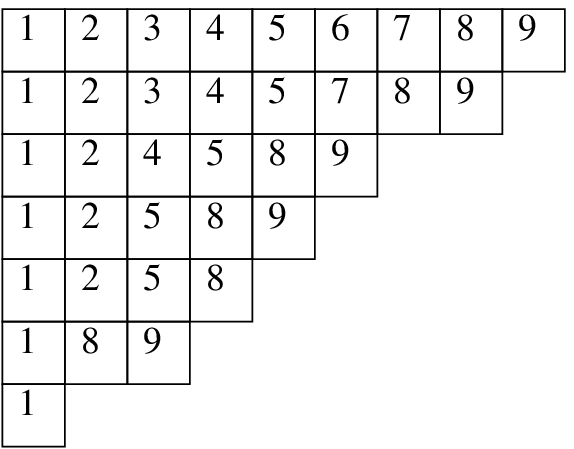} \\[.4cm]
$L_R(d)$
\includegraphics[scale=.5]{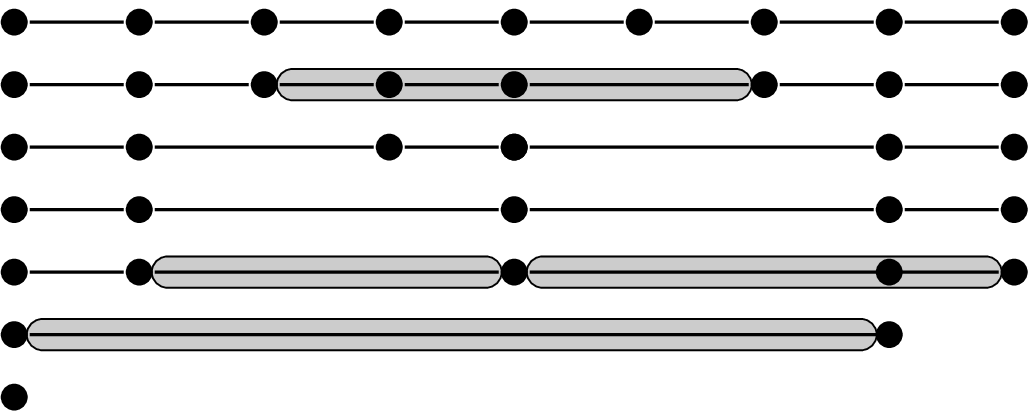} 
\end{center}
\caption{The tableaux $T(d)$, $T(i,j)$ and $L_R(d)$ for 
Example~\ref{ex:T(d)}.}\label{fig:T(d)}
\end{figure}
\end{ex}


%
%

\begin{lemma}\label{lm:T-gleich-Z}
Let $d=(d_1,\dots,d_t)$ be a dimension vector, $(i,j)\in\Gamma(d)$. 
Then 
$$
\frn(T(i,j))= Z_{ij}\,. 
$$
In particular, $Z_{i,j}$ is
irreducible.
\end{lemma}

\begin{proof}
We show that $\frn(T(i,j))=\{A\in \frn\mid\rk A[ij]^{\kappa(i,j)}< r_{ij}^{\kappa(i,j)}\}=Z_{ij}$. 
The second equation holds as it is the definition of $Z_{ij}$. \\
We first prove the lemma for a special case: replace $d_1,\ldots, d_{i-1}$ and
$d_{j+1},\ldots,d_t$ by zero, thus we get a new shorter dimension vector $e :=
(d_i,\ldots,d_j) = (e_1,\ldots,e_{j-i+1})$. Note that $(i,j)$ is in $\Gamma(d)$
precisely when $(1,j-i+1)$ is in $\Gamma(e)$. 
Also note that the codimension of
$Z_{i,j}$ for $d$ coincides with the codimension of $Z_{1,j-i+1}$ for $e$, the
first variety is just a product of the latter with an affine space.
Consequently, $Z_{i,j}$ for $d$ is irreducible precisely when $Z_{1,j-i+1}$ is
irreducible for $e$. Finally, we compare the component $\frn (T(i,j))$
for $d$ with the unique component $\frn(T(1,j-i+1))$ for $e$ that coincides
with $\mathfrak n \cap C(\mu(1,j-i+1))$ for $e$. Again, both are just given by
the equation $\rk A[i,j]^{\kappa(i,j)} < r_{i,j}^{\kappa(i,j)}$ for $d$,
respectively $\rk A[1,j-i+1]^{\kappa(1,j-i+1)} < r_{1,j-i+1}^{\kappa(1,j-i+1)}$ for $e$. This
finally shows that both varieties coincide.
%

\end{proof}

%
\section{The irreducible components of $Z$}
%

We are now ready to 
finish the proof of the descriptions of the decomposition of the 
complement 
$Z=\frn\setminus\cO(d)$ of the Richardson orbit into irreducible components. 
Again, let $d=(d_1,\dots,d_t)$ be a dimension vector, $\lambda=\lambda(d)$ 
the partition of the Richardson orbit and $(i,j)$ a pair 
with $1\le i<j\le t$. 
Recall that the $T(i,j)$ are elements of $\cT(\mu(i,j),d)$. By 
Proposition~\ref{prop:comps} the $T(i,j)$ correspond to irreducible components 
of $\frn\cap C(\mu(i,j))$. 
So the corresponding $\frn(T(i,j))$ are irreducible.

\begin{thm}\label{thm:Z-ij}
$$
Z=\bigcup_{(i,j)\in\Lambda(d)} Z_{ij}
$$
is the decomposition of $Z$ into irreducible components. 
\end{thm}

\begin{proof}
We know that $Z$ is the union of all $Z_{ij}$ over all $(i,j)$ with $1\le i<j\le t$ 
from Lemma~\ref{lm:Z-union}. 
By Lemma~\ref{lm:not-Gamma}, 
$$
Z=\bigcup_{(k,l)\in\Gamma'(d)} Z_{kl}
$$
for some subset $\Gamma'(d)\subseteq\Gamma(d)$. 
And finally, Corollary~\ref{cor:not-Lambda} tells us that for each $(k,l)$ in
this subset $\Gamma'(d)$, there exists $(i,j)\in\Lambda(d)$ such that 
$Z_{kl}$ is contained in  $Z_{ij}$. 

It remains to see that $Z_{ij}\subsetneq Z_{kl}$ and $Z_{ij}\supsetneq Z_{kl}$ for all 
$(i,j)\ne(k,l)$ $\in\Lambda(d)$. This follows as for $(i,j)\ne (k,l)$ from $\Lambda(d)$, 
one can find matrices $A$ in $Z_{ij}$ which do not satisfy the conditions 
for $Z_{kl}$ and vice versa: 
Assume $(i,j)\ne(k,l)\in \Lambda(d)$. From the line diagram $L_R(d)$ we remove 
one edge of the lowest chain connecting columns $i$ and $j$, connecting the 
resulting edges if possible with lower rows to the left and right (as with 
the dashed lines in Figure~\ref{fig:not-Gamma}) 
produces an element $A$ of $Z_{ij}$ (under $\Phi$) with $A\notin Z_{kl}$. 
It is completely analogous to find $B\in Z_{kl}$, $B\notin Z_{ij}$ 

The irreducibility follows now since $Z_{ij}=\frn(T(i,j))$ (Lemma~\ref{lm:T-gleich-Z}). 
\end{proof}

\begin{cor}\label{cor:components}
The complement $Z=\frn\setminus \cO(d)$ 
has at most $t-1$ irreducible components. 
\end{cor}

\begin{proof}
If $d$ is increasing or decreasing then clearly, $\Lambda(d)$ has size $t-1$, 
cf. Example~\ref{ex:Lambda}. 
The same is true if the $d_i$ are all different. 
In all other cases there are $d_i=d_j$ with $|j-i|>1$, and such 
that there exists an index $i<l<j$ with $d_l\ne d_i$. If $d_l>d_i$ is 
minimal among these, then neither $(i,l)$ nor $(l,j)$ belong to $\Lambda(d)$ 
and thus $\Lambda(d)$ has at most $t-2$ elements. 
The same is true for $d_l<d_i$, $d_l$ maximal among such. 
\end{proof}

Furthermore, we can describe the codimension of $Z_{ij}$ in $\frn$ as follows. Recall 
that $T(i,j)$ is obtained from $T(d)$ through a minimal movement 
(see Subsection~\ref{ss:young-tab}). 
Let $c(i,j)$ be the number of rows the box with label $j$ moves down, i.e. 
$j$ goes from row $s(i,j)$ to row $s(i,j)+c(i,j)$. 
It is known that for every row a box in a Young diagram is moved down, the 
dimension of the GL$_n$-orbit of the corresponding nilpotent 
elements decreases by two. This can be seen using the formula for the dimension of 
the stabilizer from~\cite{kp}. 
The change in dimension in the nilradical is 
half of this. Thus, 
the resulting 
$\frn(T(i,j))$ then has codimension $c(i,j)$ in the nilradical $\frn$ and we get: 

\begin{cor}\label{cor:codim}
For $(i,j)\in\Gamma(d)$, $Z_{ij}$ has codimension $c(i,j)$ in $\frn$. 
\end{cor}

The second description of the irreducible components of $Z$ is now an 
immediate consequence of Theorem~\ref{thm:Z-ij} and 
Lemma~\ref{lm:T-gleich-Z}: 

\begin{cor} \label{cor:tableaux}
$$
Z=\bigcup_{(i,j)\in\Lambda(d)} \frn(T(i,j))
$$
is the decomposition of $Z$ into irreducible components. 
\end{cor}

%
\section{An application}
%

In the last section, we illustrate our work on an example. 
We work with $G=\GL_5$ and consider the parabolic subgroups 
of different dimension vectors. \\

\noindent
A) If $d=(1,1,1,1,1)$ then $P=B$ is a Borel subgroup. Note that 
$\Lambda(d)=\Gamma(d)=\{(1,2),(2,3),(3,4),(4,5)\}$, so Theorem~\ref{thm:Z-ij}  
describes the complement $Z$ as the union 
\[
Z= Z_{12}\cup Z_{23}\cup Z_{34}\cup Z_{45}
\] 
of four irreducible components. 
 
In this example, we have that $A_{ij}=a_{ij}$ are all 
$1\times 1$-matrices. 
The Richardson 
orbit is the intersection of the regular nilpotent orbit with the set of upper 
triangular matrices in $\frgl_5$. The regular nilpotent elements are the 
nilpotent $5\times 5$-matrices whose 4th power is non-zero. 
So the Richardson orbit consists of the strictly upper triangular matrices 
$A=(a_{ij})_{ij}$ with 
\begin{eqnarray*}
A[1,5]^4 & = & \begin{pmatrix} 0 & 0 & 0 & 0 & a_{12}a_{23}a_{34}a_{45} \\ 
 & 0 & 0 & 0 & 0\\ & & 0 & 0 & 0 \\ & & & 0 & 0 \\ &&&& 0 \end{pmatrix}
  \quad \mbox{with $a_{12}a_{23}a_{34}a_{45}\ne 0$.} 
\end{eqnarray*}
For $A$ to be in the complement $Z$ of the Richardson orbit, 
the product $a_{12}a_{23}a_{34}a_{45}$ has to be zero, i.e. $A[1,5]^4=0$. 
Then clearly, $A\in Z_{i,i+1}$ 
for an $i\le 4$ as $Z_{i,i+1}$ is the set of matrices with $A_{i,i+1}=0$. 
Thus, $A$ lies in one of the components $Z_{ij}$ with $(i,j)\in \Lambda(d)$. \\

\noindent
B) If $d=(1,1,1,2)$ then $\Lambda(d)=\Gamma(d)=\{(1,2),(2,3),(3,4)\}$. 
The Richardson 
orbit is determined by the conditions 
$\rk A[12]=\rk A[23]=\rk A[34]=1$, $\rk A[13]^2=\rk[24]^2=1$, 
$\rk[14]^3=1$ (for $A\in \frn$). For $A$ to be in the complement, one 
of these ranks has to be zero. By Theorem~\ref{thm:Z-ij}, we should 
have 
$$
Z=Z_{12}\cup Z_{23}\cup Z_{34}
$$ 
where the component $Z_{12}$ consists of the matrices $A\in\frn$ with 
$a_{12}=0$, the component $Z_{23}$ of the $A$ with $a_{23}=0$ and 
$Z_{34}$ of the $A$ with $a_{34}=a_{35}=0$. 
Let us first compute $A^2$, and $A^3$ 
for $A\in\frn$ (we omit the zero entries in the opposite nilradical): 
\[
A=\begin{pmatrix}
0 & a_{12} & a_{13} & a_{14} & a_{15} \\
 & 0 & a_{23} & a_{24} & a_{25} \\
 & & 0 & a_{34} & a_{35} \\ 
 & & & 0 &0 \\
 & & & & 0
\end{pmatrix} 
\]
\[
A^2=\begin{pmatrix}
0 & 0 & a_{12}a_{23} & a_{12}a_{24}+a_{13}a_{34} & a_{12}a_{25}+a_{13}a_{35} \\
 & 0 & 0 & a_{23}a_{34} & a_{23}a_{35} \\
 & & 0 & 0 & 0 \\ 
 & & & 0 & 0 \\
 & & & & 0
\end{pmatrix}
\]

\[
A^3=\begin{pmatrix}
0 & 0 & 0 & a_{12}a_{23}a_{34}  & a_{12}a_{23}a_{35}  \\
 & 0 & 0 & 0 & 0 &  \\
 & & 0 & 0 & 0 \\ 
 & & & 0 & 0\\
 & & & & 0
\end{pmatrix}
\]
Then we see that $A[14]^3=A^3=0$ if and only if $a_{12}a_{23}a_{34}=0$ {\em and} 
$a_{12}a_{23}a_{35}=0$. Thus, $A$ clearly belongs to 
one of the three components described 
above. Now, $A[13]^2=0$ if and only if $a_{12}a_{23}=0$ as this is the only 
non-zero entry of $A[13]^2$. Similarly, $A[24]^2=0$ if and only if 
$a_{23}a_{34}=0$ {\em and} $a_{23}a_{35}=0$. 
In all cases, $A$ is contained in one of the three components.
The case of $d=(2,1,1,1)$ is completely analogous. \\

\noindent
C) The first interesting case appears for $d=(1,1,2,1)$. Here, 
$\Lambda(d)=\{(1,2),(2,4)\}\ne \Gamma(d)$. So we expect two irreducible 
components, $Z_{12}$ as the matrices $A$ with $a_{12}=0$ and 
$Z_{24}$ as the $A$ with $\rk A[24]^2=0$. We first compute $A$, $A^2$ 
and $A^3$ for $A\in \frn$: 
\[
A=\begin{pmatrix}
0 & a_{12} & a_{13} & a_{14} & a_{15} \\
 & 0 & a_{23} & a_{24} & a_{25} \\
 & & 0 & 0 & a_{35} \\ 
 & & & 0 & a_{45} \\
 & & & & 0
\end{pmatrix} \ \ \ 
A^2=\begin{pmatrix}
0 & 0 & a_{12}a_{23} & a_{12}a_{24} & a_{12}a_{25}+a_{13}a_{35}+a_{14}a_{45} \\
 & 0 & 0 & a_{23}a_{35} & a_{24}a_{45} \\
 & & 0 & 0 & 0 \\ 
 & & & 0 & 0 \\
 & & & & 0
\end{pmatrix}
\]

\[
A^3=\begin{pmatrix}
0 & 0 & 0 & 0 &a_{12}(a_{23}a_{34} +  a_{24}a_{45})  \\
 & 0 & 0 & 0 & 0 &  \\
 & & 0 & 0 & 0 \\ 
 & & & 0 & 0\\
 & & & & 0
\end{pmatrix}
\]
The elements $A$ of the Richardson orbit have non-zero  
$a_{12}$, and $\rk A[23]=\rk A[34]=1$,  
$\rk A[13]^2=\rk A[24]^2=\rk A[14]^3=\rk A^3=1$. 
Clearly, when $a_{12}=0$, then $A\in Z_{12}$. And when $\rk A[23] \rk A[34]=0$, 
$A$ belongs to $Z_{24}$. 
Now $A[14]^3=0$ if and only if $a_{12}=0$ or $a_{23}a_{34}+a_{24}a_{34}=0$ which 
is equivalent to $A\in Z_{12}$ or $A\in Z_{24}$, respectively. 
Furthermore, $A[13]^2=0$ if and only if $a_{12}a_{23}=0$ {\em and} 
$a_{12}a_{24}=0$, which is equivalent to $A\in Z_{12}\cup Z_{24}$. 
The matrices $A$ satisfying $A[24]^2$ are by definition $Z_{24}$. 
The case $d=(1,2,1,1)$ is analogous. \\

\noindent
D) Let $d=(2,2,1)$, with $\Lambda(d)=\Gamma(d)=\{(1,2),(2,3)\}$, the complement 
should be $Z_{12}\cup Z_{23}$. 
The Richardson 
orbit is given as the matrices $A$ with $\rk A[12]=2$ and $\rk A[23]=\rk A[13]^2=1$. 
We compute $A$ and $A^2$: 
\[
A=\begin{pmatrix}
0 &  0 & a_{13} & a_{14} & a_{15} \\
 & 0 & a_{23} & a_{24} & a_{25} \\
 & & 0 & 0 & a_{35} \\ 
 & & & 0 & a_{45} \\
 & & & & 0
\end{pmatrix} \ \ \ 
A^2=\begin{pmatrix}
0 & 0 & 0 & 0 & a_{13}a_{35}+a_{14}a_{45} \\
 & 0 & 0 & 0 & a_{23}a_{35} + a_{24}a_{45} \\
 & & 0 & 0 & 0 \\ 
 & & & 0 & 0 \\
 & & & & 0
\end{pmatrix}
\]
If $A$ is a matrix with $A[13]^2=0$ then if $a_{35}=a_{45}=0$, 
$A$ is an element of $Z_{23}$. So let $\rk A[23]\ne 0$. 
Solving the two equations $a_{13}a_{35}+a_{14}a_{45}=0$ 
$a_{23}a_{35} + a_{24}a_{45}=0$ then shows that the rank 
of $A[12]$ is one.  Thus $Z_{13}$ is already contained in $Z_{12}$. 
The case $d=(1,2,2)$ is analogous.\\

\noindent
E) The second interesting case is $d=(2,1,2)$, with 
$\Lambda(d)=\{(1,3)\}$ and $\Gamma(d)=\{(1,3), (1,2), (2,3)\}$. 
Here we only obtain one irreducible component 
in the complement! 
The Richardson orbit is defined by $\rk A[13]^2=\rk A^2=1$ 
and $\rk A=3$: The dimension of its stabilizer has to be equal to 
the dimension of the Levi factor. Using the formulae from~\cite{kp} then gives 
this description of the Richardson orbit. 
For the complement, we are looking at matrices $A$ with 
$\rk A[12]=0$ or $\rk A[23]=0$ or $\rk A[13]^2=0$. 
If $A$ satisfies $A[12]=0$ then $A^2$ is also zero, so 
$A\in Z_{13}$ by definition. Similarly, matrices with $A[23]=0$ square to zero 
and hence lie in $Z_{13}$. \\

\noindent
F) 
The case $d=(1,3,1)$ with $\Lambda(d)=\{(1,3)\}$, so again, we only
have one component in the complement of the open dense orbit. 
For matrices of the Richardson orbit, we have $\rk A[12]=\rk A[23]=\rk A[13]^2=1$. 
For the complement, we take matrices where one of these ranks is 
zero. 
If it is $\rk A[12]=0$ or $\rk A[23]=0$ then clearly, $A[13]^2=0$, so $A\in Z_{13}$. 
The cases $d=(3,1,1)$ and $d=(1,1,3)$ behave similarly as $d=(2,1,1,1)$ 
and $d=(1,1,1,2)$. We omit them here. \\

\noindent
G) 
The remaining cases are $d=(4,1)$, $d=(1,4)$. Here, the complement to 
the Richardson orbit is given by $A[12]=0$, i.e. it is the zero matrix.


%

\end{document}